\documentclass[12pt,a4paper]{amsart}
\usepackage{amsmath,amssymb,amscd}
\input xy
\xyoption{all}

\newcommand{\ra}{\rightarrow}

\newcommand{\RR}{\mathbb R}
\newcommand{\CC}{\mathbb C}
\newcommand{\ZZ}{\mathbb Z}

\newcommand{\cH}{\mathcal{H}}

\newcommand{\cO}{\mathcal{O}}
\newcommand{\cM}{\mathcal{M}}
\newcommand{\IM}{\mbox{Im}}
\newcommand{\NS}{\mbox{NS}}
\newcommand{\Pic}{\mbox{Pic}}
\newcommand{\End}{\mbox{End}}
\newcommand{\Aut}{\mbox{Aut}}
\newcommand{\Hom}{\mbox{Hom}}
\newcommand{\diag}{\mbox{diag}}
\theoremstyle{plain}
\newtheorem{theorem}{Theorem}[section]
\newtheorem{lemma}[theorem]{Lemma}
\newtheorem{proposition}[theorem]{Proposition}
\newtheorem{cor}[theorem]{Corollary}

\begin{document}
\title[Products of elliptic curves]{Principal polarizations on products of elliptic curves}

\author{Herbert Lange}

\address{H. Lange\\Mathematisches Institut\\
              Universit\"at Erlangen-N\"urnberg\\
              Bismarckstra\ss e $1\frac{ 1}{2}$\\
              D-$91054$ Erlangen\\
              Germany}
             
\email{lange@mi.uni-erlangen.de}
                  \dedicatory{Dedicated to Sevin Recillas}

\thanks{Supported by DFG Contract Ba 423/8-1}

\begin{abstract}
An abelian variety admits only a finite number of isomorphism classes of principal polarizations. The paper gives an 
interpretation of this number in terms of class numbers of definite Hermitian forms in the case of a product of elliptic curves 
without complex multiplication. In the case of a self-product of an elliptic curve, as well as in the two-dimensional case, 
classical class number computations can be applied to determine this number.  
\end{abstract}
\maketitle

\section{Introduction}

\noindent
Let $X$ be an abelian variety of dimension $g$ over the field $\CC$ of complex numbers. 
The {\it N\'eron-Severi group} of $X$ is by definition the quotient of the group of line bundles on $X$
modulo the subgroup of line bundle which are algebraically equivalent to zero,
$$
\NS(X) = \Pic(X)/\Pic^0(X).
$$ 
A {\it polarization} of $X$ is by definition an element $l \in \NS(X)$ which is represented by an ample line bundle
$L$ on $X$. By a slight abuse of notation we denote the polarization defined by $L$ often by $L$ itself. The pair $(X,L)$ is called a {\it polarized 
abelian variety}. Every polarization $L$ of $X$ induces an 
isogeny  $\varphi_L: X \ra \hat{X}$ of $X$ onto its dual abelian variety $\hat{X} = \Pic^0(X)$, defined by
$$
\varphi_L(x) = t^*_xL \otimes L^{-1},
$$  
where $t_x: X \ra X$ denotes the translation by $x$. A polarization on $X$ is called {\it principal}, if 
$\varphi_L$ is an isomorphism. 

Two polarizations $L_1$ and $L_2$ on $X$ are said to be {\it equivalent} if there is an automorphism 
$\tau \in \mbox{Aut(X)}$ such that
$$
\tau^*L_2 \sim L_1,
$$
where  $\sim$  denotes algebraic equivalence. Consider the set
$$
P(X) = \{ \emph{equivalence classes of principal polarizations on}\, X \}.
$$ 
According to a theorem of Narasimhan-Nori \cite{NN},
$P(X)$ is a finite set. It is the aim of this note to compute the number $\# P(X)$ in the case of a 
product of elliptic curves without complex multiplication,
$$
X = E_1 \times \cdots \times E_n.
$$ 
According to Lemma \ref{lemma2.2} we may assume that the $E_i$ are pairwise isogenous. The main result is Theorem \ref{thm3.5},
which says that in this case there is a bijection between $P(X)$ and the set of equivalence classes of Hermitian forms 
of a certain type. 

We can then apply results on class numbers of Hermitian or quadratic forms in order to compute the cardinality of
the set $P(X)$ in certain cases. If for example $E$ is an elliptic curve without complex multiplication and
$$
X = E \times \cdots \times E,
$$
then
$$
\# P(X) = h(n)
$$
(see Theorem \ref{thm4.1}), where $h(n)$ denotes the number of equivalence classes of positive
definite integral quadratic forms of 
rank $n = \dim X$. For $n \leq 25$,
the class number $h(n)$ is known. We deduce that for $n \leq 7$ the abelian variety $X$ admits no principal polarization apart from
the canonical one. For $n=8$, there is another principal polarization and we have for example $\# P(X) = 8$ for $n=16$ 
and $\# P(X) = 297$ for $n=24$.
Moreover, it is a consequence of the mass formula of Minkowski-Siegel that the number $\# P(X)$ is unbounded for $n \ra \infty$.

For another application of Theorem \ref{thm3.5} consider
$$
X = E_1 \times E_2
$$ 
with isogenous elliptic curves $E_1$ and $E_2$ without complex multiplication admitting an isogeny of minimal positive degree $d$. Here
there is a bijection between $P(X)$ and the set of equivalence classes of primitive positive definite integral quadratic forms
of rank 2 and determinant $d$.
The corresponding class number $\widetilde{h}(d)$ has been computed by Hayashida in \cite{ha}, also in order to compute the number 
of classes of principal polarizations of $E_1 \times E_2$ as above, but with a different approach. It is a consequence of this result 
that the number
$\# P(X)$ is unbounded for $d \ra \infty$.\\

The main idea for the proof of Theorem \ref{thm3.5} is as follows: The canonical principal polarization of $X$ induces a 
bijection between $P(X)$ and the set of equivalence classes of symmetric automorphisms of $X$. Via the analytic representation of 
$X$ and a suitable choice of bases this set can be considered as a set of equivalence classes of Hermitian forms.  

\section{Generalities}

\noindent
Let $X$ denote a complex abelian variety of dimension $n$. In order to compute the number $\# P(X)$ we may assume that 
$X$ admits at least one principal polarization, say $L_0$, which is fixed in the sequel. 
This polarization induces an anti-involution on the ring $\End(X)$ of endomorphisms of $X$, the Rosati involution defined by
$$
\alpha \mapsto \alpha' = \varphi_{L_0}^{-1} \hat{\alpha} \varphi_{L_0}, 
$$
for every $\alpha \in \End(X)$, where $\hat{\alpha}: \hat{X} \ra \hat{X}$ denotes the dual endomorphism. Denote by
$$
\End^s(X) = \{ \alpha \in \End(X) \;|\; \alpha' = \alpha \}
$$
the abelian subgroup of symmetric endomorphisms. Similarly $\Aut^s(X)$ is defined. Finally, set 
$$
\Aut^s(X)^+ = \{ \alpha \in \Aut^s(X) \;|\; \alpha \; \mbox{totally positive}\},
$$
where totally positive means that all roots of the minimal polynomial are positive.  The group $\Aut(X)$ acts on 
the set $\Aut^s(X)^+$ by
$$
(\varphi, \alpha) \mapsto \varphi' \alpha \varphi
$$
for all $\varphi \in \Aut(X)$ and $\alpha \in \Aut^s(X)^+$.
For the proof of the following proposition we refer to \cite{cav}, 
Proposition 5.2.1 and Theorem 5.2.4.

\begin{proposition} \label{prop2.1}
With the notation above we have
\begin{itemize}
\item[(1)] The map 
$$
\epsilon: \NS(X) \ra \End^s(X), \;\; L \mapsto \varphi_{L_0}^{-1} \varphi_L
$$
is an isomorphism of groups.
\item[(2)] The map $\epsilon$ induces a bijection
$$
P(X) \stackrel{\sim}{\ra} \Aut^s(X)^+/_{\sim} \;, 
$$  
where $\Aut^s(X)^+/_{\sim}$ denotes the set of equivalence classes with respect to the above action.
\end{itemize} 
\end{proposition}
\noindent
A polarization $L$ on $X$ is called {\it reducible} if $L = L_1 \otimes L_2$ with ample line bundles $L_1$ and $L_2$, or 
equivalently, if there are abelian subvarieties $X_1$ and $X_2$ of $X$ such that there is an isomorphism of 
polarized abelian varieties 
$$
(X,L) \simeq (X_1,L_1) \times (X_2,L_2),
$$
where here $L_i$ denotes the restriction of $L$ to $X_i$ for $i=1, 2$.

\begin{lemma} \label{lemma2.2}
Let $(X_i,L_i)$ be principally polarized abelian varieties for $i=1,2$. If $\Hom(X_1,X_2) = 0$, then every principal 
polarization of $X_1 \times X_2$ is reducible.
\end{lemma}

\begin{proof}
The assumption implies $\End(X) = \End(X_1) \oplus \End(X_2)$. The principal polarization $p_1^*L_1 \otimes p_2^*L_2$
induces an isomorphism $\NS(X) \stackrel{\sim}{\ra} \End^s(X) = \End^s(X_1) \oplus \End^s(X_2)$. Together with  
Proposition \ref{prop2.1} (2) this gives the assertion.
\end{proof}

Now let $E_i$ be complex elliptic curves and consider the abelian variety
$$
X = E_1 \times \cdots \times E_n.
$$
It admits the canonical principal polarization 
$$
L_0 = p_1^* \cO_{E_1}(0) \otimes \cdots \otimes p_n^* \cO_{E_n}(0),
$$
where $p_i : X \ra E_i$ denotes the $i$-th projection. Lemma \ref{lemma2.2} implies that, in order to determine all principal 
polarizations of $X$, we may assume that the elliptic curves $E_i$ are pairwise isogenous.

For the computation of $\#P(X)$ we introduce suitable period matrices of $X$.    
For $i=1, \ldots, n$ there is an element $z_i$ in the upper half plane $\cH$ such that
$$
E_i = \CC/\Lambda_i \quad \mbox{with} \quad \Lambda_i = \ZZ + z_i \ZZ.
$$
Then the following matrix is a period matrix for the abelian variety $X$ with respect to the canonical basis $e_1, \ldots, e_n$ of
$\CC^n$ and a suitable basis $\lambda_1, \ldots, \lambda_{2n}$ of the lattice $\Lambda$:
$$
\Pi = ( I_n \;Z),
$$
where $I_n$ denotes the unit matrix of degree $n$ and 
$$
Z = \diag(z_1, \ldots, z_n).
$$
satisfies
\begin{equation} \label{eq1}
Z^t = Z \quad \mbox{and} \quad \IM \,Z > 0.
\end{equation}

The first Chern class of the polarization $L_0$ can be considered as an alternating form $E_{L_0}$ on the lattice $\Lambda$.
With respect to the basis $\lambda_1, \cdots,\lambda_{2n}$ 
of $\Lambda$ the alternating form $E_{L_0}$ is given by the 
following matrix,
also denoted by $E_{L_0}$,
$$
E_{L_0} = \left( \begin{array}{cc}
                 0 & I_n\\
                 -I_n& 0
                 \end{array}  \right).  
$$
For any $\varphi \in \End(X)$ the analytic and rational representations
$$
A_{\varphi} := \rho_a(\varphi) \in M_n(\CC) \quad \mbox{and} \quad R_{\varphi} := \rho_r(\varphi) \in M_{2n}(\ZZ)
$$
satisfy the following equation
\begin{equation} \label{eq2}
A_{\varphi} \Pi = \Pi R_{\varphi}. 
\end{equation}
Conversely, any pair of matrices $(A,R) \in M_n(\CC) \times M_{2n}(\ZZ)$ satisfying (\ref{eq2}) defines an endomorphism of $X$.
The following lemma computes the analytic representation of the Rosati involution of $\varphi$ in terms of the chosen bases.

\begin{lemma} \label{lemma2.3}
\hspace{1cm} $A_{\varphi'} = \IM \,Z \cdot \overline{A_{\varphi}}^t \cdot (\IM \,Z)^{-1}.$
\end{lemma} 

\begin{proof}
It is well-known (see \cite{cav}, Proposition 5.1.1) that the rational representation $R_{\varphi'}$ of $\varphi'$ is the adjoint 
matrix of $R_{\varphi}$ with respect to the alternating form $E_{L_0}$, which implies
\begin{equation} \label{eq3}
R_{\varphi'} = E_{L_0}^{-1} R^t_{\varphi} E_{L_0}
\end{equation}
The matrix $\left( \begin{array}{c}
                    \Pi\\
                    \overline{\Pi}
                    \end{array}   \right)$
is invertible, since $\Pi$ is a period matrix. Hence (\ref{eq2}) applied to $\varphi$ and $\varphi'$ implies
\begin{equation} \label{eq4}
R_{\varphi} = \left( \begin{array}{c}
                    \Pi\\
                    \overline{\Pi}
                    \end{array}   \right)^{-1}
\left( \begin{array}{cc}
        A_{\varphi}&0\\
        0& \overline{A_{\varphi}}
        \end{array} \right)
        \left( \begin{array}{c}
                    \Pi\\
                    \overline{\Pi}
                    \end{array}   \right)
\end{equation}
and
\begin{equation} \label{eq5}
\left( \begin{array}{cc}
        A_{\varphi'}&0\\
        0& \overline{A_{\varphi'}}
        \end{array} \right) = \left( \begin{array}{c}
                    \Pi\\
                    \overline{\Pi}
                    \end{array}   \right)  R_{\varphi'}
                                                       \left( \begin{array}{c}
                    \Pi\\
                    \overline{\Pi}
                    \end{array}   \right)^{-1}
\end{equation} 
Inserting (\ref{eq3}) and the transpose of (\ref{eq4}) into (\ref{eq5}), we obtain
$$
\left( \begin{array}{cc}
        A_{\varphi'}&0\\
       0 & \overline{A_{\varphi'}}
        \end{array} \right) = S \left( \begin{array}{cc}
        A_{\varphi}^t&0\\
        0& \overline{A_{\varphi}^t}
        \end{array} \right) S^{-1}
$$        
with 
$$
S= \left( \begin{array}{c}
                    \Pi\\
                    \overline{\Pi}
                    \end{array}   \right)
                    \left( \begin{array}{cc}
                    0& -I_n\\
                    I_n&0
                    \end{array} \right)
                    \left( \begin{array}{cc}
                    \Pi^t&
                    \overline{\Pi^t}
                    \end{array}   \right) = 2i \left( \begin{array}{cc}
                                                      0& \IM \,Z\\
                                                      - \IM \,Z&0
                                                      \end{array} \right).
$$
For the last equation we used (\ref{eq1}). This implies the assertion.
\end{proof}

\section{$E_i$ without complex multiplication}

\noindent
Let the notation be as at the end of the last section. In particular $E_i = \CC/\Lambda, \; i=1, \ldots, n$ are 
pairwise isogenous complex elliptic curves. In this section we assume in addition that the $E_i$ are 
without complex multiplication which means that $\End(E) \simeq \ZZ$.  Then 
$$
\Hom(E_i,E_j) \simeq \ZZ
$$
for all $i,j = 1, \cdots , n$. Denote by $\tau_{ij}: E_i \ra E_j$ an isogeny of minimal positive degree, say $d_{ij}$. 
We identify $\tau_{ij}$ with its analytic 
representation, i.e., we consider it as a complex number. 

\begin{lemma} \label{lem3.1}
The complex conjugate $\overline{\tau_{ij}}$ represents a homomorphism $E_j \ra E_i$, denoted by the same symbol, such that
$$
\overline{\tau_{ij}} \tau_{ij} = \deg({\tau_{ij}})1_{E_i} \quad \mbox{and} \quad 
\tau_{ij} \overline{\tau_{ij}} = \deg({\tau_{ij}})1_{E_j}.
$$
$\overline{\tau_{ij}}$ is a homomorphism of minimal positive degree.
\end{lemma}

\begin{proof}
There exists a homomorphism 
$\widetilde{\tau}_{ji}: E_j \ra E_i$ such that 
$\widetilde{\tau_{ji}} \tau_{ij} = \deg({\tau_{ij}})1_{E_i}$ and 
$\tau_{ij} \widetilde{\tau_{ji}} = \deg({\tau_{ij}})1_{E_j}$.
In particular
$$
\widetilde{\tau_{ji}} \tau_{ij} = \deg(\tau_{ij}).
$$
 On the other hand, 
since the rational representation is the direct sum of the analytic representation and its complex conjugate,
$$
\overline{\tau_{ij}} \tau_{ij} = \det(R_{\tau_{ij}}) = \deg(\tau_{ij}).
$$
This implies that the analytic representation of the homomorphism $\widetilde{\tau_{ij}}$ is given by the complex number $\overline{\tau_{ij}}$
and thus completes the proof of the lemma.
\end{proof} 

Any homorphism $\varphi_{ij}: E_i \ra E_j$ can be written
as
$$
\varphi_{ij} = d_{ij} \tau_{ij}
$$
with a uniquely determined integer $d_{ij}$. This remark together with Lemma \ref{lem3.1} imply the following proposition.

\begin{proposition} 
The analytic representation induces an isomorphism 
$$
\End(X) \stackrel{\simeq}{\ra} \cM_n(X), \quad \varphi \mapsto A_{\varphi}
$$
where $\cM_n(X)$ denotes the ring of all matrices of the form
\begin{equation}
A_{\varphi} = \begin{pmatrix}
              d_{11}&d_{12}\overline{\tau_{12}}& \cdots& d_{1n}\overline{\tau_{1n}}\\
              d_{21}\tau_{12}&d_{22}& \cdots&d_{2n}\overline{\tau_{2n}}\\
              \vdots&\vdots&\ddots&\vdots\\
              d_{n1}\tau_{1n}&d_{n2}\tau_{2n}&\cdots&d_{nn}\\
              \end{pmatrix}
\end{equation}
with $\tau_{ij}$ as above and $d_{ij} \in \ZZ$ for all $i,j$. 
\end{proposition}

\begin{lemma} \label{lemma3.2} For any $\varphi \in \End(X)$ we have, with respect to the chosen bases,

$$
\det A_{\varphi} \in \ZZ.
$$
\end{lemma}

\begin{proof}
The matrix is of the form $A_{\varphi} = (\varphi_{ij})$ with homomorphisms $\varphi_{ij}: E_i \ra E_j$, which are 
identified with their analytic representations. 
By definition of the determinant,
$$
\det A_{\varphi} = \sum_{\sigma \in S_n} sign(\sigma) \varphi_{1\sigma(1)}\varphi_{2\sigma(2)} \cdots \varphi_{n\sigma(n)}.
$$
Now any $\sigma \in S_n$ is a product of cycles $\sigma = \sigma_r \cdots \sigma_1$. But any cycle $\sigma_{\nu}$ 
of length $s_{\nu}$ say, is of the form  
$$
\sigma_{\nu} = (k,\sigma(k),\sigma^2(k),\cdots, \sigma^{s_{\nu} - 1}(k)).
$$ 
Hence the corresponding homomorphism $\varphi_{k\sigma(k)}\varphi_{\sigma(k)\sigma^2(k)} \cdots \varphi_{\sigma^{s_{\nu}-1}(k),k}$
 is an endomorphism of $E_k$. So by assumption any cycle is an element of $\End(E_i) = \ZZ$ for some $i$. 
This implies that $\det A_{\varphi}$ is an integer as a sum of products of integers. 
\end{proof}

As in the last section, $\Pi = (I_n,Z)$ with $Z= \diag(z_1,\cdots ,z_n)$
is a period matrix for $X = E_1 \times \cdots \times E_n$.
If $\{e_1, \ldots, e_n\}$ denotes the canonical basis of $\CC^n$, we introduce a new basis $\{f_1, \ldots, f_n\}$ of 
$\CC^n$ by setting 
$$
f_i = \sqrt{\IM z_i}e_i
$$
for $i=1,\ldots,n$. Since $\IM(z_i) > 0$ for all $i$, 
$$
T = \diag(\sqrt{\IM(z_1)}, \ldots ,\sqrt{\IM(z_n)})
$$
is a well-defined matrix of $GL_n(\RR)$. The period matrix of $X$ with respect to the new basis $\{f_i\}$ of $\CC^n$ and 
the old basis $\{\lambda_j\}$ of $\Lambda$ is  
$$
\widetilde{\Pi} = T^{-1}\Pi.
$$
Using this it is an easy consequence of Lemma \ref{lemma2.3} that the analytic representation $A_{\varphi'}$ of the
Rosati transform $\varphi'$ of $\varphi$ with respect to these bases is given by
$$
A_{\varphi'} = \overline{A_{\varphi}}^t.
$$
Thus we can conclude

\begin{lemma} \label{lemma3.4}
For $\varphi \in End(X)$ the following conditions are equivalent
\begin{itemize}
\item[(1)] $\varphi$ is symmetric with respect to the Rosati involution,
\item[(2)] The analytic representation $A_{\varphi}$ with respect to the bases $\{f_i\}$ of $\CC^n$ and $\{\lambda_j\}$
of $\Lambda$ is a Hermitian matrix.
\end{itemize}  
\end{lemma}

Note that Lemma \ref{lemma3.4} does not assume that the elliptic curves are without complex multiplication. For the proof of the lemma
a normalization, different from the one used above, turns out to be more convenient (see Remark 3.6).\\ 

Consider now the set
$$
\cM_n^s(X)^+ = \{ A \in \cM_n(X) \;|\; \overline{A}^t = A, \; A > 0,\; \det A = 1 \}.
$$
Two matrices $A_1, A_2 \in \cM_n^s(X)^+$ are called {\it equivalent}, if there is an invertible matrix $T \in \cM_n(X)$ such that
$$
A_2 = \overline{T}^t A_1 T.
$$ 
Note that this equivalence is just the matrix version of the usual equivalence of Hermitian forms. Recall that $P(X)$ denotes the set of 
equivalence classes of principal polarizations of $X$.

\begin{theorem} \label{thm3.5}
There is a bijection 
$$
P(X) \ra \cM_n^s(X)^+/_{\sim},
$$
where $\cM_n^s(X)^+/_{\sim}$ denotes the set of equivalence classes of Hermitian matrices in $\cM_n^s(X)^+$.
\end{theorem}

\begin{proof}
We claim first that the analytic representation with respect to the basis $\{f_i\}$ induces a bijection 
$Aut^s(X)^+ \ra \cM^s_n(X)^+$. 

For the proof note that by Lemma \ref{lemma3.4} an endomorphism $\alpha$ of $X$ is symmetric if and only if $A_{\alpha}$ is Hermitian. It is
totally positive if and only if all zeros of the minimal polynomial are positive, i.e., $A_{\alpha}$ is positive definite.
On the other hand, $\alpha$ is an automorphism if and only if $\deg(\alpha) = 1$. But $\deg(\alpha)$ equals the determinant of the 
rational representation $R_{\alpha}$. So (\ref{eq4}) implies 
$$
\det A_{\alpha} \cdot \overline{\det A_{\alpha}} = 1.
$$
By Lemma \ref{lemma3.2}, $\det A_{\alpha} \in \ZZ$, which implies $\det A_{\alpha} = \pm1$ and thus $= 1$, if
$A_{\alpha}$ is positive definite. This implies the assertion.

It is clear that the bijection $Aut^s(X)^+ \ra \cM^s_n(X)^+$ is compatible with the equivalence relations. Hence 
Proposition \ref{prop2.1} completes the proof of the theorem.
\end{proof}

\noindent
{\bf Remark 3.6.} One can use the same method in order to prove an analogous result for a product of pairwise 
isogenous elliptic curves with complex multiplication. This has been done in the thesis of P. Schuster,
written under the supervision of the author (see \cite{sch}).
Again the idea is to use suitable bases in order to interpret the set $P(X)$ in terms of classes of Hermitian forms. 
In the case of dimension 2 and of self-products of dimension 3 Schuster applies class number formulas of Hashimoto-Koseki
in order to compute the number $\#P(X)$.

\section{Self-products of an elliptic curve $E$}

\noindent
Let $E$ be an elliptic curve over $\CC$ without complex multiplication and consider the abelian variety
$$
X = E \times \cdots \times E
$$
of dimension $n$. In this case we may choose $\tau_{ij} = 1_E$ for all $i,j$. Then the following theorem 
is a special case of Theorem \ref{thm3.5}.

\begin{theorem} \label{thm4.1}
There is a bijection of $P(X)$ with the set of equivalence classes of positive definite unimodular integral 
quadratic forms. 
\end{theorem}

\begin{cor} \label{cor3.7}
For $n \leq 7$, there is no principal polarization on $X$ apart from $L_0$.
\end{cor}

\begin{proof}
For $n \leq 7$, there is only one positive definite unimodular integral quadratic form (see \cite{se}).
\end{proof}

There is an extensive literature about the number $h(n)$ of classes of positive definite unimodular integral quadratic forms of rank $n$
(see \cite{CS} and the literature quoted there).
According to Theorem \ref{thm4.1} this number can be interpreted as the number $\# P(X)$ of equivalence classes of principal 
polarizations of $X$. In particular $h(n)$ has been computed for $n \leq 25$ by Kneser, Niemeier, Conway-Sloane and Borchards 
(see \cite{CS}, table 2.2).
Together with Theorem \ref{thm4.1} this gives  

\begin{cor} 
\ \hfill
\begin{tabular}[t]{r||r|r|r|r|r|r|r|r|r|} \hline
n & 8 & 9 & 10 & 11 & 12 & 13 & 14 & 15 & 16\\ \hline
$\#P(X)$ & 2 & 2 & 2 & 2 & 3 & 3 & 4 & 5 & 8\\ \hline \hline
n & 17 & 18 & 19 & 20 & 21 & 22 & 23 & 24 & 25\\ \hline
$\#P(X)$ & 9 & 13 & 16 & 28 & 40 & 68 & 117 & 297 & 665\\ \hline
\end{tabular}
\end{cor}

The mass formula of Siegel-Minkowski gives an estimate for the number $h(n)$ (see \cite{CS} or \cite{se}). 
As a consequence one gets for example that $h(32) \geq 80.000.000$ and moreover
that $h(n)$
tends to infinity if $n \ra \infty$. \\

\noindent
{\bf Remark 4.4.} The theory of unimodular definite integral quadratic forms distinguishes between even and odd ones 
(or of type I and type II).
It would be interesting to see whether this distinction has a geometric meaning for the corresponding principal polarizations.

\section{Abelian surfaces $E_1 \times E_2$}

\noindent
In this section consider 
$$
X = E_1 \times E_2
$$
with isogenous elliptic curves $E_1$ and $E_2$ without complex multiplication and assume that 
$$
\tau: E_1 \ra E_2
$$ 
is an isogeny of minimal 
positive degree $d \geq 2$.
According to Theorem \ref{thm3.5} there is a bijection between $P(X)$ and the set of classes of Hermitian forms 
$\cM_2^s(X)^+/_{\sim}$, where
$$
\cM_2^s(X)^+ = \{ A = \begin{pmatrix}
                   a_{11}&a_{12}\overline{\tau}\\
                   a_{12}\tau& a_{22}\\
                   \end{pmatrix}  \; | \; a_{ij} \in \ZZ, \; a_{11} > 0, \; \det A = 1 \}.
$$ 

In order to compute the corresponding class number, recall from Lemma \ref{lem3.1} that the complex conjugate 
$\overline{\tau}$ represents an isogeny $E_2 \ra E_1$ of minimal positive degree $d$, and consider the map
$$
\Phi: End(E_1 \times E_2) \ra End(E_1 \times E_1),\;\; \;\;\varphi \mapsto (1_{E_1} \times \overline{\tau}) \varphi (1_{E_1} \times \tau)
$$ 
In terms of the analytic representation with respect to the basis $\{f_j\}$ of section 3, the map $\Phi$ is given by
$$
A_{\varphi} \mapsto \Phi(A_{\varphi}) = \begin{pmatrix} 1&0\\ 0 & \overline{\tau} \end{pmatrix}
         \begin{pmatrix} a_{11}& a_{12}\overline{\tau}\\ a_{12} \tau& a_{22} \end{pmatrix}
                                        \begin{pmatrix} 1&0\\ 0 & \tau \end{pmatrix}
       =  \begin{pmatrix} a_{11}& a_{12}d\\ a_{12} d& a_{22}d \end{pmatrix}.                                          
$$
So $\Phi(A_{\varphi})$ is an integral quadratic form of determinant $d$. 
Recall that a integral $(2 \times 2)$-matrix $(m_{ij})$ is called {\it primitive}, if $\gcd(m_{ij}\; |\; i,j=1,2) = 1$.    
Consider the set of integral quadratic forms given as matrices by
$$
\widetilde{\cM}_2^s(X)^+ = \{ B = \begin{pmatrix} b_{11}& b_{12}\\ b_{12}& b_{22} \end{pmatrix} \;|\; 
                                   B \; \mbox{primitive}, \; b_{11} >0, \; \det{B} = d \}.
$$
As usual, two integral quadratic forms 
$B_1$ and $B_2$ are called {\it equivalent} if there is a $T \in GL_2(\ZZ)$ such that $B_2 = T^tB_1T$. 
Clearly this defines an equivalence relation on the set $\widetilde{\cM}_2^s(X)^+$. It is well-known
 that there are only finitely many equivalence classes of such forms (in fact, this is 
also a consequence of the following theorem). Let $\tilde{h}(d)$ denote the corresponding class number:
$$
\tilde{h}(d) = \#(\widetilde{\cM}_2^s(X)^+/\sim).
$$

\begin{theorem} \label{thm3.8}
Let $X = E_1 \times E_2$ with elliptic curves $E_1$ and $E_2$ admiiting an isogeny of minimal positive degree $d$. Then
$$
\quad \# P(X) = \tilde{h}(d).
$$
\end{theorem}

\begin{proof}
According to Theorem \ref{thm3.5} it is sufficient to show that the map $\Phi:\cM_2^s(X)^+ \ra \widetilde{\cM}_2^s(X)^+$
is compatible with the equivalence relations and induces a bijection of the equivalence classes.

Let $A = \begin{pmatrix} a_{11} & a_{12}\overline{\tau}\\
                         a_{12}\tau& a_{22} \end{pmatrix} \in \cM_2^s(X)^+$. Then
$\Phi(A) = \begin{pmatrix} a_{11} & a_{12}d\\
                         a_{12}d& a_{22}d \end{pmatrix}$ is a primitive matrix, since
$$
\gcd(a_{11},a_{12}d,a_{22}d) = \gcd(a_{11},\gcd(a_{11},a_{22})d) = \gcd(a_{11},d) = 1,
$$
where we used the fact that $\det(A) = 1$.

Next we claim that for $A_1, A_2 \in \cM_2^s(X)^+$,
$$
A_1 \sim A_2 \quad \Leftrightarrow \quad \Phi(A_1) \sim \Phi(A_2).
$$
Let $T = \begin{pmatrix} t_{11}& t_{12}\overline{\tau}\\ t_{21}\tau& t_{22} \end{pmatrix}$ be an invertible matrix
of the ring $\cM_2(X)$ with $A_2 = T^tA_1T$. Then $\Phi(A_2) = \widetilde{T}^t \Phi(A_1) \widetilde{T}$ with 
$\widetilde{T} = \begin{pmatrix} t_{11}& t_{12}d\\t_{21}&t_{22} \end{pmatrix}$. This implies $\Phi(A_1) \sim \Phi(A_2)$, since
$\det \widetilde{T} = \det T$.

Conversely, given $B_1 = \begin{pmatrix} b_{11}& b_{12}d\\b_{12}d&b_{22}d \end{pmatrix} = \Phi(A_1)$ and 
$B_2 = \begin{pmatrix} c_{11}& c_{12}d\\ c_{12}d& c_{22}d \end{pmatrix} \newline = \Phi(A_2)$ with $B_2 = T^tB_1T, \; 
T = \begin{pmatrix}t_{11}&t_{12}\\t_{21}&t_{22} \end{pmatrix}$. This means in particular
$$
c_{11} = t_{11}^2b_{11} + d(2t_{11}t_{21}b_{12} + t_{21}^2b_{22})
$$
$$
c_{12}d = t_{11}t_{12}b_{11} + d(t_{12}t_{21}b_{12} + t_{11}t_{22}b_{12} + t_{21}t_{22}b_{22}) 
$$
Since $B_1$ and $B_2$ are primitive, this implies
$$
t_{12} \equiv 0 \mod d,
$$ 
from which we get $A_1 \sim A_2$ by reading the above computation upside down.

Hence $\Phi$ induces a map $\overline{\Phi}:\cM_2^s(X)^+/_\sim \ra \widetilde{\cM}_2^s(X)^+/_\sim$. It remains to show that 
$\overline{\Phi}$ is bijective. 

Clearly $\overline{\Phi}$ is injective. To see that it is  surjective, let 
$B = (b_{ij}) \in \widetilde{\cM}_2^s(X)^+$. According to an elementary result for binary quadratic forms (see e.g., 
\cite{M}, p. 132) there are integers $x_0, y_0$, prime to each other, such that
$\gcd(b_{11}x_0^2 + 2b_{12}x_0y_0 + b_{22}y_0^2,d) =1$. Choose $x_1,y_1 \in \ZZ$ with $x_0y_1-y_0x_1 =1$. Using this we see that,
replacing $B$ by $U^tBU$ with $U=\begin{pmatrix} x_0&x_1\\y_0&y_1 \end{pmatrix}$, we may assume $gcd(b_{11},d) = 1$.
So $b_{11}r + ds=1$ with integers $r,s$. Let $T=\begin{pmatrix} 1&-b_{12}r\\0&1 \end{pmatrix}$. Then
$$
B \sim T^tBT = \begin{pmatrix} b_{11}& b_{12}sd\\b_{12}sd & b_{11}b_{12}^2r^2 - 2b_{12}^2r + b_{22} \end{pmatrix} \in \IM(\Phi),
$$
since $d$ divides $b_{12}$ and $b_{22}$. This means that $\overline{\Phi}$ is surjective.
\end{proof}

The class number $\tilde{h}(d)$ has been computed by T. Hayashida in \cite{ha}, actually for the same purpose as here, namely the 
determination of the number of classes of principal polarization on $E_1 \times E_2$ as above. However the 
way to associate a quadratic form of $\widetilde{\cM}_2^s(X)^+$ to a principal polarization and thus the proof of 
Theorem \ref{thm3.8} was completely different. I will not repeat the actual class numbers here, but refer to \cite{ha} for details.
Note only that the formulas imply that the number $\#P(X)$ is unbounded as $d \ra \infty$.

\end{document}